\documentclass[a4paper,11pt]{amsart}
\usepackage{amsmath,amsthm,amssymb}
\usepackage[mathscr]{eucal}
 \usepackage{cite}
\usepackage{upgreek}
\usepackage[bookmarks,bookmarksnumbered]{hyperref}
\usepackage{enumerate}
\usepackage{amsmath}
\usepackage{mathrsfs}
\usepackage{blindtext}
\usepackage{scrextend}
\usepackage{enumitem}
\addtokomafont{labelinglabel}{\sffamily}
\usepackage{color}

\numberwithin{equation}{section}

\theoremstyle{plain}

\newtheorem{theorem}{Theorem}[section]

\newtheorem{proposition}[theorem]{Proposition}
\newtheorem{corollary}[theorem]{Corollary}
\theoremstyle{definition}
\newtheorem{definition}[theorem]{Definition}
\newtheorem*{definition*}{Definition}
\newtheorem{example}[theorem]{Example}

\newtheorem{remark}[theorem]{Remark}

\begin{document}

\title[Rigidity for von Neumann algebras]
{Rigidity for von Neumann algebras}

\author[A. Ioana]{Adrian Ioana}
\address{Department of Mathematics, University of California San Diego, 9500 Gilman Drive, La Jolla, CA 92093, USA, and IMAR, Bucharest, Romania}
\email{aioana@ucsd.edu}
\thanks{A.I. was supported in part by NSF Career Grant DMS \#1253402 and a Sloan Foundation Fellowship.}
\begin{abstract} 

We survey some of the progress made recently in the classification of von Neumann algebras arising from countable groups and their measure preserving actions on probability spaces.  
We emphasize results which provide classes of (W$^*$-superrigid) actions that can be completely recovered from their von Neumann algebras and II$_1$ factors that have a unique Cartan subalgebra. We also present cocycle superrigidity theorems and some of their applications to orbit equivalence. Finally, we discuss several recent rigidity results for von Neumann algebras associated to groups.

\end{abstract}

\maketitle



\maketitle

\section{Introduction}
\noindent
\subsection{}
A {\bf von Neumann algebra} is an algebra of bounded linear operators on a Hilbert space which is closed under the adjoint operation and in the weak operator topology.  Von Neumann algebras arise naturally from countable groups and their actions on probability spaces, via two seminal constructions of Murray and von Neumann \cite{MvN36,MvN43}. Given a countable group $\Gamma$,  the left regular representation of $\Gamma$ on $\ell^2\Gamma$ generates the {\bf group von Neumann algebra} $L(\Gamma)$. Equivalently, $L(\Gamma)$  is the weak operator closure of the complex group algebra $\mathbb C\Gamma$ acting on $\ell^2\Gamma$ by left convolution. 
Every nonsingular action $\Gamma\curvearrowright (X,\mu)$ of a countable group $\Gamma$ on a probability space $(X,\mu)$ gives rise to the {\bf group measure space von Neumann algebra} $L^{\infty}(X)\rtimes\Gamma$.

A central theme in the theory of von Neumann algebras is the classification of $L(\Gamma)$ in terms of the group $\Gamma$ and of $L^{\infty}(X)\rtimes\Gamma$  in terms of the group action $\Gamma\curvearrowright (X,\mu)$. 
These problems are typically studied when $\Gamma$ has infinite non-trivial conjugacy classes (icc) and when $\Gamma\curvearrowright (X,\mu)$ is free ergodic and measure preserving, respectively. These assumptions guarantee that the corresponding algebras are {\bf II$_1$ factors}: indecomposable infinite dimensional von Neumann algebras which admit a  trace. Moreover, it follows that $L^{\infty}(X)$ is a {\bf Cartan subalgebra} of $L^{\infty}(X)\rtimes\Gamma$, that is, a maximal abelian subalgebra whose normalizer generates $L^{\infty}(X)\rtimes\Gamma$.

The classification of II$_1$ factors is a rich subject,
 with deep connections  to several areas of mathematics.  
 Over the years, it has repeatedly provided fertile ground for the development of new, exciting theories: Jones' subfactor theory, Voiculescu's free probability theory, and Popa's deformation/rigidity theory. 
 The subject has been connected to group theory and ergodic theory from its very beginning, via the group measure space construction \cite{MvN36}.  Later on, Singer made the 
 observation that the isomorphism class of $L^{\infty}(X)\rtimes\Gamma$ only depends on the equivalence relation given by the orbits of $\Gamma\curvearrowright (X,\mu)$ \cite{Si55}. This soon led to a new branch of ergodic theory, which studies group actions up to orbit equivalence \cite{Dy58}. 
 Orbit equivalence theory, further developed in the 1980s (see \cite{OW80,CFW81,Zi84}), has seen an explosion of activity in the last twenty years (see  \cite{Sh04,Fu09,Ga10}).  This progress has been in part triggered by the success of the deformation/rigidity approach to the classification of II$_1$ factors (see \cite{Po06b,Va10a,Io12a}). More broadly, this approach has generated a wide range of applications to ergodic theory and descriptive set theory, including:

\begin{itemize}
\item the existence of non-orbit equivalent actions of non-amenable groups \cite{GP03,Io06b,Ep07}.
\vskip 0.05in
\item  cocycle superrigidity theorems for Bernoulli actions in \cite{Po05,Po06a}, leading to examples of non-Bernoullian factors of Bernoulli actions for a class of countable groups \cite{Po04b,PS03}, and the solution of some open problems in descriptive set theory \cite{Th07}
\vskip 0.05in
\item cocycle superrigidity theorems for profinite actions \cite{Io08, Fu09, GITD16}, leading to an explicit uncountable family of Borel incomparable treeable equivalence relations \cite{Io13}.
\vskip 0.05in
\item solid ergodicity of Bernoulli actions \cite{CI08} (see also \cite{Oz04}).
\end{itemize}

\subsection{}
The classification of II$_1$ factors is governed by a strong amenable/non-amenable dichotomy.
Early work in this area culminated with Connes' celebrated theorem from the mid 1970s: all II$_1$ factors arising from infinite amenable groups and their actions are isomorphic to the hyperfinite II$_1$ factor \cite{Co75b}. Amenable groups thus manifest a striking absence of rigidity: any property of the group or action, other than the amenability of the group, is lost in the passage to von Neumann algebras.

In contrast, it gradually became clear that in the non-amenable case various aspects of groups and actions are remembered by their von Neumann algebras. Thus, non-amenable groups were used in \cite{Mc69} and \cite{Co75a} to construct large families of non-isomorphic II$_1$ factors. Rigidity phenomena for von Neumann algebras  first emerged in the work of Connes from the early 1980s \cite{Co80}. He showed that II$_1$ factors arising from property (T) groups have countable symmetry (fundamental and outer automorphism) groups. 
Representation theoretic properties of groups (Kazhdan's property (T), Haagerup's property, weak amenability) were then used to prove unexpected non-embeddability results for II$_1$ factors associated to certain lattices in Lie groups  \cite{CJ83, CH88}.
But while these results showcased the richness of the theory, the classification problem for non-amenable II$_1$ factors remained by and large intractable.

A major breakthrough in the classification of II$_1$ factors was made by Popa with his invention of deformation/rigidity theory \cite{Po01b,Po03,Po04a} (see the surveys \cite{Po06b,Va06a}). The studied II$_1$ factors, $M$, admit a distinguished subalgebra $A$ (e.g., $L^{\infty}(X)$ or $L(\Gamma)$ when $M=L^{\infty}(X)\rtimes\Gamma$) such that the inclusion $A\subset M$ satisfies both a deformation and a rigidity property. Popa discovered that the combination of these properties can  be used to detect the position of $A$ inside $M$, 
or even recover the underlying structure of $M$ (e.g., the group $\Gamma$ and action $\Gamma\curvearrowright X$ when $M=L^{\infty}(X)\rtimes\Gamma$). 
He also developed a series of powerful technical tools to exploit this principle.

Popa first implemented this idea and techniques to provide a class of II$_1$ factors, $M$, which admit a unique Cartan subalgebra subalgebra, $A$, with the relative property (T) \cite{Po01b}.  
The uniqueness of $A$ implies that any invariant of the inclusion $A\subset M$ is in fact  an invariant of $M$.
When applied to  $M=L^{\infty}(\mathbb T^2)\rtimes$ SL$_2(\mathbb Z)$, it follows that the fundamental group of $M$ is equal to that of the equivalence relation of the action SL$_2(\mathbb Z)\curvearrowright \mathbb T^2$.
Since the latter is trivial by Gaboriau's work \cite{Ga99,Ga01}, this makes $M$ the first example of a II$_1$ factor with trivial fundamental group \cite{Po01b}, thereby solving a longstanding problem. 

 In \cite{Po03,Po04a}, Popa greatly broadened the scope of deformation/rigidity theory by obtaining the first strong rigidity theorem for group measure space factors.  
To make this precise, suppose that $\Gamma\curvearrowright (X,\mu)$ is a Bernoulli action of an icc group $\Gamma$ and $\Lambda\curvearrowright (Y,\nu)$ is a free ergodic probability measure preserving action of a property (T) group $\Lambda$ (e.g.,  $\Lambda=\text{SL}_{n\geq 3}(\mathbb Z)$).
 Under this assumptions, it is shown in \cite{Po03,Po04a} that if the group measure space factors $L^{\infty}(X)\rtimes\Gamma$ and $L^{\infty}(Y)\rtimes\Lambda$ are isomorphic, then the groups $\Gamma$ and $\Lambda$ are isomorphic and their actions are conjugate.

\subsection{}
The goal of this survey is to present some of the progress achieved recently in the classification of II$_1$ factors. We focus on advances made since 2010, and refer the reader to \cite{Po06b,Va10a} for earlier developments.  A topic covered there but omitted here is the calculation of symmetry groups of II$_1$ factors; see \cite{Po01b,Po03, IPP05, PV08a} for several key results in this direction.  There is some overlap with \cite{Io12a}, but overall the selection of topics and presentation are quite different. 

We start with a section of preliminaries (Section 2) and continue with a discussion of the main ideas from Popa's deformation/rigidity theory (Section 3).  
Before giving an overview of Sections 4-6, we first recall some terminology and background. 
Two free ergodic p.m.p. (probability measure preserving) actions $\Gamma\curvearrowright (X,\mu)$ and $\Lambda\curvearrowright (Y,\nu)$ are called:
\begin{enumerate}
\item {\bf conjugate} if there exist an isomorphism of probability spaces $\alpha:(X,\mu)\rightarrow (Y,\nu)$ and an isomorphism of groups $\delta:\Gamma\rightarrow\Lambda$ such that $\alpha(g\cdot x)=\delta(g)\cdot\alpha(x)$, for all $g\in\Gamma$ and almost every $x\in X$.
\vskip 0.05in
\item {\bf orbit equivalent (OE)} if there exists an isomorphism of probability spaces \\ $\alpha:(X,\mu)\rightarrow (Y,\nu)$ such that $\alpha(\Gamma\cdot x)=\Lambda\cdot\alpha(x)$, for almost every $x\in X$.
\vskip 0.05in
\item {\bf W$^*$-equivalent} if $L^{\infty}(X)\rtimes\Gamma$ is isomorphic to $L^{\infty}(Y)\rtimes\Lambda$.
\end{enumerate}
Singer showed that OE amounts to the existence of an isomorphism  $L^{\infty}(X)\rtimes\Gamma\cong L^{\infty}(Y)\rtimes\Lambda$ which identifies the Cartan subalgebras $L^{\infty}(X)$ and $L^{\infty}(Y)$ \cite{Si55}. Thus, OE implies W$^*$-equivalence.  
Since conjugacy clearly implies orbit equivalence, putting these together we have:
 $$\text{conjugacy}\;\;\Longrightarrow\;\;\text{orbit equivalence}\;\;\Longrightarrow\;\;\text{W$^*$-equivalence}$$

Rigidity usually refers to a situation in which a weak equivalence between two objects
can be used to show that the objects are equivalent in a much stronger sense or even
isomorphic. In the present context, rigidity occurs whenever some of the above implications can be reversed for all actions $\Gamma\curvearrowright (X,\mu)$ and $\Lambda\curvearrowright (Y,\nu)$ belonging to two classes of actions. 
The most extreme form of rigidity, called superrigidity, happens when this can be achieved without any restrictions on the  second class of actions. 
Thus, an action $\Gamma\curvearrowright (X,\mu)$ is  {\bf W$^*$-superrigid} (respectively, {\bf OE-superrigid}) if any free ergodic p.m.p. action $\Lambda\curvearrowright (Y,\nu)$ which is W$^*$-equivalent (respectively, OE) to  $\Gamma\curvearrowright (X,\mu)$ must be conjugate to it. In other words, the conjugacy class of the action can be entirely reconstructed from the isomorphism class of its von Neumann algebra (respectively, its orbit equivalence class).

The seminal results from \cite{Po01b,Po03,Po04a} suggested two rigidity conjectures which have guided much of the work in the area in the ensuing years.
First, \cite{Po01b} provided a class of II$_1$ factors $L^{\infty}(X)\rtimes\mathbb F_n$ associated to actions of the free groups $\mathbb F_n$, with $n\geq 2$, for which $L^{\infty}(X)$ is the unique Cartan subalgebra  satisfying the relative property (T). This led to the conjecture that the same might be true for arbitrary Cartan subalgebras of arbitrary free group measure space factors: 
{\bf (A)} {\it $L^{\infty}(X)\rtimes\mathbb F_n$  has a unique Cartan subalgebra, up to unitary conjugacy,  for  any  free ergodic p.m.p. action $\mathbb F_n\curvearrowright (X,\mu)$ of $\mathbb F_n$ with $n\geq 2$.}
Second,  \cite{Po03,Po04a} showed that within the class of Bernoulli actions of icc property (T) groups, W$^*$-equivalence implies conjugacy.  Moreover, it was proved in \cite{Po05} that such actions are OE-superrigid. These results naturally led to the following conjecture: 
{\bf (B)} {\it Bernoulli actions $\Gamma\curvearrowright (X,\mu)$ of icc property (T) groups $\Gamma$ are W$^*$-superrigid}.

In Section 4, we discuss the positive resolutions of the above conjectures. These were the culmination of a period of intense activity which has generated a series of striking unique Cartan decomposition and W$^*$-superrigidity results.
The first breakthrough was made by Ozawa and Popa who confirmed conjecture {\bf (A)} in the case of profinite actions  \cite{OP07}.
The class of groups whose profinite actions give rise to II$_1$ factors with a unique group measure space decomposition was then shown to be much larger in \cite{OP07,OP08,Pe09}.  

But, since none of these actions was known to be OE-superrigid, W$^*$-superrigidity could not be deduced. 
Indeed, proving that an action $\Gamma\curvearrowright (X,\mu)$ is W$^*$-superrigid amounts to showing that the action is OE-superrigid {\it and} that $L^{\infty}(X)\rtimes\Gamma$ has a unique {\it group measure space Cartan subalgebra}\footnote{This terminology is used to distinguish the Cartan subalgebras coming from the group measure space construction from general Cartan subalgebras, see Section \ref{cartan}. }.   
Nevertheless, Peterson was able to show  the existence of ``virtually" W$^*$-superrigid actions \cite{Pe09}.
Soon after, Popa and Vaes discovered a large class of amalgamated free product groups $\Gamma$ whose every free ergodic p.m.p. action $\Gamma\curvearrowright (X,\mu)$ gives rise to a II$_1$ factor with a unique group measure space Cartan subalgebra, up to unitary conjugacy \cite{PV09}.
 Applying OE-superrigidity results from \cite{Po05,Po06a,Ki09} enabled them to provide the first concrete classes of W$^*$-superrigid actions.
However, these results do not apply to actions of property (T) groups and so despite all of this progress, conjecture {\bf (B)} remained open until it was eventually settled by the author in \cite{Io10}. 

New uniqueness theorems for group measure space Cartan subalgebras were then obtained in \cite{CP10} by Chifan and Peterson and in \cite{Io11a}, while Chifan and Sinclair extended the results of \cite{OP07}  from free groups to hyperbolic groups in \cite{CS11}.  
However, the available uniqueness results for Cartan subalgebras  required either a rigidity property of the group (excluding the free groups), or the action to be in a specific class. 
Thus, the situation for arbitrary actions of the free groups remained unclear until Conjecture {\bf (A)} was resolved by Popa and Vaes  in their breakthrough work \cite{PV11}. Subsequently,  several additional families of groups were shown to satisfy  conjecture {\bf (A)}, including non-elementary hyperbolic groups \cite{PV12}, free product groups \cite{Io12a}, and central quotients of braid groups \cite{CIK13}, while conjecture {\bf (B)} was established for the more general class of mixing Gaussian actions in \cite{Bo12}.

Section 5 is devoted to rigidity results in orbit equivalence.  
After recalling the pioneering rigidity and superrigidity phenomena discovered by Zimmer \cite{Zi80,Zi84} and Furman \cite{Fu98}, we discuss several recent cocycle superrigidity results.  Our starting point is Popa's remarkable cocycle superrigidity theorem: any cocycle for a Bernoulli action of a property (T) group  into a countable group is cohomologous to a homomorphism \cite{Po05}. 
The property (T) assumption was removed in \cite{Po06a}, where the same was shown to hold for Bernoulli actions of products of non-amenable groups. 
A cocycle superrigidity theorem for a different class of actions of property (T) groups, the profinite actions, was then obtained in \cite{Io08}. Answering a question motivated by the analogy with Bernoulli actions, this theorem was recently extended to product groups in \cite{GITD16}. 

These and many other results provide large classes of ``rigid" groups (including property (T) groups, product groups, and by Kida's work \cite{Ki06}, most mapping class groups) which admit OE-superrigid actions. In contrast, other non-amenable groups,  notably the free groups $\mathbb F_n$, posses no OE-superigid actions.
Nevertheless, as we explain in the second part of Section 4, a general OE-rigidity theorem for profinite actions was discovered in \cite{Io13}.  This result imposes no assumptions on the acting groups and so it applies, in novel fashion, to actions of groups such as $\mathbb F_n$ and SL$_2(\mathbb Z)$. 
As an application, it led to a continuum of mutually non-OE and Borel incomparable actions of SL$_2(\mathbb Z)$, confirming a conjecture from \cite{Th01,Th06}. It also motivated a ``local spectral gap" theorem for dense subgroups of simple Lie groups in  \cite{BISG15}. 

In Section 6, we discuss recent rigidity results for group von Neumann algebras. These give instances when certain algebraic properties of groups, such as the absence or presence of a direct product decomposition, are  remembered by their von Neumann algebras. A remarkable result of Ozawa shows that for any icc hyperbolic group $\Gamma$, the II$_1$ factor $L(\Gamma)$ is prime, i.e. it cannot be decomposed as a tensor product of two II$_1$ factors \cite{Oz03}. In particular, this recovered the primeness of $L(\mathbb F_n)$, for $n\geq 2$, which was first proved in \cite{Ge96} using Voiculescu's free probability techniques. Later on, several other
 large classes of 
  icc groups $\Gamma$ were shown to give rise to prime II$_1$ factors (see e.g. \cite{Oz04, Pe06, Po06c,CH08}).  
 However, all such groups, $\Gamma$, 
satisfy various properties which relate them closely to lattices in rank one simple Lie groups. On the other hand, it is an  open problem whether II$_1$ factors associated to icc irreducible lattices $\Gamma$ in higher rank simple or semisimple Lie groups $G$ are prime. We present in Section 6 a result from
\cite{DHI16} which answers this positively in the case when $G$ is a product of simple Lie groups of rank one.

We continue with the recent finding  in \cite{CdSS15,CI17} of large classes of product and amalgamated free product groups whose product (respectively, amalgam) structure can be recognized from their von Neumann algebras. 
Finally, we turn our attention to the strongest type of rigidity for group II$_1$ factors $L(\Gamma)$, called W$^*$-superrigidity. This occurs when the isomorphism class of $\Gamma$ can be  reconstructed from the isomorphism class of $L(\Gamma)$. 
 We conclude the section with a discussion of the first examples of W$^*$-superrigid groups discovered in our joint work with Popa and Vaes \cite{IPV10}, and the subsequent examples exhibited in \cite{BV12,CI17}.

\subsection{} 
Let us mention a few exciting topics related to the classification of II$_1$ factors which have received a lot of attention recently but are not covered here, due to limitations of space. First, 
we point out the impressive work of Houdayer and his co-authors
  (including \cite{HV12,BHR12, HI15, BHV15}) where the deformation/rigidity framework is adapted to study von Neumann algebras of type III.
A notable advance in this direction is the classification theorem for free Araki-Woods factors obtained by Houdayer, Shlyakhtenko, and Vaes in \cite{HSV16}. 
 We also highlight Peterson's remarkable work  \cite{Pe15} (see also \cite{CP13})  which shows that lattices in higher rank simple Lie groups admit a unique II$_1$ factor representation, the regular representation,  thus solving a conjecture of Connes from the 1980s. 
 Finally, we mention the model theory for II$_1$ factors which was introduced in \cite{FHS09,FHS10,FHS11}. Subsequently, our  joint work with Boutonnet and Chifan \cite{BCI15} settled a basic question in the theory. More precisely, we showed the existence of uncountably many different elementary classes of II$_1$ factors (equivalently, of uncountably many II$_1$ factors with pairwise non-isomorphic ultrapowers). 

\subsection{Acknowledgements} It is my pleasure to thank R\'{e}mi Boutonnet, Cyril Houdayer, and Sorin Popa for many comments that helped improve the exposition.


\section{Preliminaries}

\subsection{Tracial von Neumann algebras} 
A von Neumann algebra $M$ is called {\it tracial} if it admits a linear functional $\uptau:M\rightarrow\mathbb C$, called a {\it trace}, which is
\begin{enumerate}
\item {\it positive}:  $\uptau(x^*x)\geq 0$, for all $x\in M$.
\item {\it faithful}: $\uptau(x^*x)=0$, for some $x\in M$, implies that $x=0$.
\item {\it normal}:  $\uptau(\sum_{i\in I}p_i)=\sum_{i\in I}\uptau(p_i)$, for any family $\{p_i\}_{i\in I}$ of mutually orthogonal projections. 
\item {\it tracial}:  $\uptau(xy)=\uptau(yx)$, for all $x,y\in M$.
\end{enumerate}

A von Neumann algebra with trivial center is called a {\it factor}. An infinite dimensional tracial factor is called a {\it II$_1$ factor}. 
Note that any II$_1$ factor  $M$ admits a unique trace $\uptau$ such that $\uptau(1)=1$.

Any tracial von Neumann algebra $(M,\uptau)$ admits a  canonical (or {\it standard}) representation on a Hilbert space. Denote by $L^2(M)$  the completion of $M$ with respect to the $2$-norm $\|x\|_2:=\sqrt{\uptau(x^*x)}$. 
Then the left  and right multiplication actions of $M$ on itself give rise to representations of $M$ and its opposite algebra $M^{\text{op}}$ on $L^2(M)$. This makes $L^2(M)$ a {\it Hilbert $M$-bimodule}, i.e. a Hilbert space $\mathcal H$ endowed with commuting normal representations $M\subset\mathbb B(\mathcal H)$ and $M^{\text{op}}\subset\mathbb B(\mathcal H)$. 

Let $P\subset M$ be a von Neumann subalgebra and denote by $E_P:M\rightarrow P$ the unique $\uptau$-preserving conditional expectation onto $P$. Specifically, $E_P$ is determined by the identity $\uptau(E_P(x)y)=\uptau(xy)$, for all $x\in M$ and $y\in P$.
By completing the algebraic tensor product $M\otimes_{\text{alg}}M$ with respect to the scalar product $\langle x_1\otimes x_2,y_1\otimes y_2\rangle=\uptau(y_2^*E_P(x_2^*x_1)y_1)$ we obtain the Hilbert $M$-bimodule $L^2(M)\bar{\otimes}_PL^2(M)$. Alternatively, this bimodule can be realized as the L$^2$-space of Jones' basic construction $\langle M,P\rangle$. In the case $P=M$ and $P=\mathbb C1$, we recover the so-called {\it trivial} and {\it coarse} bimodules, $L^2(M)$ and $L^2(M)\bar{\otimes}L^2(M)$, respectively.

\subsection{Group von Neumann algebras}
Let $\Gamma$ be a countable group and denote by $\{\delta_h\}_{h\in\Gamma}$ the usual orthonormal basis of $\ell^2\Gamma$. The {\it left regular representation} $u:\Gamma\rightarrow\mathcal U(\ell^2\Gamma)$ is given by $u_g(\delta_h)=\delta_{gh}$. The {\it group von Neumann algebra} $L(\Gamma)$ is defined as the weak operator closure of the linear span of $\{u_g\}_{g\in\Gamma}$. It is a tracial von Neumann algebra with a trace $\uptau:L(\Gamma)\rightarrow\mathbb C$ given by $\uptau(x)=\langle x\delta_e,\delta_e\rangle$. Equivalently, $\uptau$ is the unique trace on $L(\Gamma)$ satisfying $\uptau(u_g)=\delta_{g,e}$, for all $g\in\Gamma$. 

Note that $L(\Gamma)$ is a II$_1$ factor if and only if $\Gamma$ is {\it icc}: $\{ghg^{-1}|g\in\Gamma\}$ is infinite, for all $h\in\Gamma\setminus\{e\}$.

\subsection{Group measure space von Neumann algebras} Let $\Gamma\curvearrowright (X,\mu)$ be a p.m.p. action of a countable group $\Gamma$ on a probability space $(X,\mu)$.
For $g\in\Gamma$ and $c\in L^{2}(X)$, let $\sigma_g(c)\in L^{2}(X)$ be defined by $\sigma_g(c)(x)=c(g^{-1}\cdot x)$.   The elements of both $\Gamma$ and $L^{\infty}(X)$ can be represented as operators on the Hilbert space $L^2(X)\bar{\otimes}\ell^2\Gamma$ through the formulae $u_g(c\otimes\delta_h)=\sigma_g(c)\otimes\delta_{gh}$ and $b(c\otimes\delta_h)=bc\otimes\delta_h$. Then $u_g$ is a unitary operator and $u_gbu_g^*=\sigma_g(b)$, for all $g\in\Gamma$ and $b\in L^{\infty}(X)$. As a consequence, the linear span of $\{bu_g|b\in L^{\infty}(X),g\in\Gamma\}$ is a $*$-algebra.

The {\it group measure space von Neumann algebra} $L^{\infty}(X)\rtimes\Gamma$ is defined as the weak operator closure of the linear span of $\{bu_g|b\in L^{\infty}(X),g\in\Gamma\}$.  
It is a tracial von Neumann algebra with a trace $\uptau:L^{\infty}(X)\rtimes\Gamma\rightarrow\mathbb C$ given by $\uptau(a)=\langle a(1\otimes\delta_e),1\otimes\delta_e\rangle$.
Equivalently, $\uptau$ is the unique trace on $L^{\infty}(X)\rtimes\Gamma$ satisfying $\uptau(bu_g)=\delta_{g,e}\int_Xb\;\text{d}\mu$, for all $g\in\Gamma$ and $b\in L^{\infty}(X)$. 

Any element $a\in L^{\infty}(X)\rtimes\Gamma$ admits a {\it Fourier decomposition} $a=\sum_{g\in\Gamma} a_gu_g$, with $a_g\in L^{\infty}(X)$ and the series converging in $\|.\|_2$. These so-called {\it Fourier coefficients} $\{a_g\}_{g\in\Gamma}$ of $a$  are given by  $a_g=E_{L^{\infty}(X)}(au_g^*)$ and satisfy $$\sum_{g\in\Gamma}\|a_g\|_2^2=\|a\|_2^2.$$ 
 
 If the action $\Gamma\curvearrowright (X,\mu)$ is (essentially) free and ergodic, then $L^{\infty}(X)\rtimes\Gamma$ is a II$_1$ factor.
 Recall that the action $\Gamma\curvearrowright (X,\mu)$ is called {\it free} if $\{x|\;g\cdot x=x\}$ is a null set, for all $g\in\Gamma\setminus\{e\}$, and {\it ergodic} if any $\Gamma$-invariant measurable subset $A\subset X$ must satisfy $\mu(A)\in\{0,1\}$.

\subsection{Examples of free ergodic p.m.p. actions} 
\begin{enumerate}
\item Let $\Gamma$ be a countable group and $(X_0,\mu_0)$ be a non-trivial probability space. Then $\Gamma$ acts on the space $X_0^{\Gamma}$ of sequences $x=(x_h)_{h\in\Gamma}$ by shifting the indices: $g\cdot x=(x_{g^{-1}h})_{h\in\Gamma}$. This action, called the {\it Bernoulli action} with base $(X_0,\mu_0)$, preserves the product probability measure $\mu_0^{\Gamma}$, and is free and ergodic.
\vskip 0.1in
\item Let $\Gamma$ be a countable group together with a dense embedding into a compact group $G$. 
The {\it left translation action} of $\Gamma$ on $G$ given by $g\cdot x=gx$ preserves the Haar measure ${\bf m}_G$ of $G$, and is free and ergodic. 
For instance, assume that $\Gamma$ is residually finite, and let $G=\varprojlim\Gamma/\Gamma_n$ 
 be its {\it profinite competition} with respect to a chain $\Gamma=\Gamma_0>\Gamma_1>...>\Gamma_n>...$ of finite index normal subgroups with trivial intersection. 
 Then $G$ is a profinite hence totally disconnected compact group, and the map $g\mapsto (g\Gamma_n)_n$ gives a dense embedding of $\Gamma$ into $G$. At the opposite end, we have the case when $G$ is a connected compact group (e.g. $G=SO(n)$) and $\Gamma$ is a countable dense subgroup of $G$.  
 \vskip 0.1in
\item Generalizing example (2), a p.m.p. action $\Gamma\curvearrowright (X,\mu)$ is called {\it compact} if the closure of $\Gamma$ in Aut$(X,\mu)$ is compact, and {\it profinite} if it is an inverse limit of actions $\Gamma\curvearrowright (X_n,\mu_n)$, with $X_n$ a finite set, for all $n$. Any profinite p.m.p. action is compact. Any ergodic compact p.m.p. action is isomorphic to a left translation action of the form $\Gamma\curvearrowright (G/K, {\bf m}_{G/K})$, where $G$ is a compact group containing $\Gamma$ densely, $K<G$ is a closed subgroup and ${\bf m}_{G/K}$ is the unique $G$-invariant Borel probability measure on $G/K$. Any ergodic profinite p.m.p. action is of this form, with $G$ a profinite group.
 \vskip 0.1in
\item Finally, the standard action of SL$_n(\mathbb Z)$ on the $n$-torus $\mathbb T^n=\mathbb R^n/\mathbb Z^n$ preserves the Lebesgue measure, and is free and ergodic. The left multiplication action of PSL$_n(\mathbb Z)$ on SL$_n(\mathbb R)/\text{SL}_n(\mathbb Z)$ preserves the unique SL$_n(\mathbb R)$-invariant probability measure, and is free and ergodic. 
 \end{enumerate}
\subsection{Equivalence relations} An equivalence relation $\mathcal R$ on a standard probability space $(X,\mu)$ is called {\it countable p.m.p.}  if $\mathcal R$ has countable classes, $\mathcal R$ is a Borel subset of $X\times X$, and any Borel  automorphism of $X$ whose graph is contained in $\mathcal R$ preserves $\mu$ \cite{FM75}.
If $\Gamma\curvearrowright (X,\mu)$ is a p.m.p. action of a countable group $\Gamma$, its {\it orbit equivalence relation} $\mathcal R_{\Gamma\curvearrowright X}:=\{(x,y)\in X^2\;|\;\Gamma\cdot x=\Gamma\cdot y\}$ is a countable p.m.p. equivalence relation. 
 It is now clear that two p.m.p. actions of countable groups $\Gamma\curvearrowright (X,\mu)$ and $\Lambda\curvearrowright (Y,\nu)$ are orbit equivalent precisely when their OE relations are isomorphic: $(\alpha\times\alpha)(\mathcal R_{\Gamma\curvearrowright X})=\mathcal R_{\Lambda\curvearrowright Y}$, for some isomorphism of probability spaces $\alpha:(X,\mu)\rightarrow (Y,\nu)$.
 
\subsection{Cartan subalgebras}\label{cartan}
Let $M$ be a II$_1$ factor. The {\it normalizer} of a subalgebra $A\subset M$, denoted by $\mathcal N_M(A)$,  is the group of unitaries $u\in M$ satisfying $uAu^*=A$.
An abelian von Neumann subalgebra $A\subset M$ is called a {\it Cartan subalgebra} if it is maximal abelian and its normalizer generates $M$. 
For instance, if $\Gamma\curvearrowright (X,\mu)$ is a free ergodic p.m.p. action, then $A:=L^{\infty}(X)$ is a Cartan subalgebra of the group measure space II$_1$ factor $M:=L^{\infty}(X)\rtimes\Gamma$. To distinguish such Cartan subalgebras from arbitrary ones we call them of {\it group measure space} type. 

In general, any Cartan subalgebra inclusion  $A\subset M$ can be identified with an inclusion of the form $L^{\infty}(X)\subset L(\mathcal R,w)$, where $(X,\mu)$ is a probability space and $L(\mathcal R,w)$  is the von Neumann algebra associated to a countable  p.m.p. equivalence relation $\mathcal R$ on $X$ and a $2$-cocycle $w\in$ H$^2(\mathcal R,\mathbb T)$ \cite{FM75}.
By \cite{FM75}, $\mathcal R$ arises as the OE relation of a p.m.p action $\Gamma\curvearrowright (X,\mu)$.  If the action is free and the cocycle is trivial, then we canonically have $M=L^{\infty}(X)\rtimes\Gamma$. However, the action cannot always be chosen to be free \cite{Fu98}, and  thus not all Cartan subalgebras are of group measure space type.

The next proposition makes clear the importance of Cartan subalgebras in the study of group measure space factors.

\begin{proposition}[Singer, \cite{Si55}]\label{singer} If $\Gamma\curvearrowright (X,\mu)$ and $\Lambda\curvearrowright (Y,\nu)$ are free ergodic p.m.p. actions, then the following conditions are equivalent
\begin{enumerate}
\item the actions $\Gamma\curvearrowright (X,\mu)$ and $\Lambda\curvearrowright (Y,\nu)$ are orbit equivalent.
\item there exists a $*$-isomorphism $\theta:L^{\infty}(X)\rtimes\Gamma\rightarrow L^{\infty}(Y)\rtimes\Lambda$ such that $\theta(L^{\infty}(X))=L^{\infty}(Y)$.

\end{enumerate}
\end{proposition}

Proposition \ref{singer}  is extremely useful in two ways. First, the implication (1) $\Rightarrow$ (2) provides an approach to the study of orbit equivalence of actions using von Neumann algebras. This has been instrumental in several developments, including the finding of non-OE actions of non-amenable groups (see \cite{GP03,Io06a, Io06b}) and of new OE-superrigidity results (see \cite{Po05,Po06a,Io08}). On the other hand, the implication (2) $\Rightarrow$ (1) allows one to reduce the classification of group measure space factors to the classification of the corresponding actions up to orbit equivalence, whenever uniqueness of group measure space Cartan subalgebras can be established. This has been used for instance to exhibit the first families of W$^*$-superrigid actions in \cite{Pe09,PV09,Io10}.

\subsection{Amenability and property (T)} 
In the early 1980s, Connes discovered that 
Hilbert bimodules provide an appropriate representation theory for tracial von Neumann algebras, paralleling the theory of unitary representations for groups (see \cite{Co82,Po86}). 

To illustrate this point, assume that $M=L(\Gamma)$, for a countable group $\Gamma$. 
Given a unitary representation $\pi:\Gamma\rightarrow\mathcal U(\mathcal H)$ on a Hilbert space $\mathcal H$,  the Hilbert space $\mathcal H\bar{\otimes}\ell^2\Gamma$ carries a natural Hilbert $M$-bimodule structure: $u_g(\xi\otimes\delta_h)u_k=\pi(g)(\xi)\otimes\delta_{ghk}$. 
Moreover, the map $\xi\mapsto \xi\otimes\delta_e$ turns sequences of $\Gamma$-almost invariant unit vectors into sequences of $M$-almost central  tracial vectors. 
Here, for a Hilbert $M$-bimodule $\mathcal K$ and a subalgebra $P\subset M$, we say that a vector $\xi\in\mathcal K$ is {\it tracial} if $\langle x\xi,\xi\rangle=\langle \xi x,\xi\rangle=\uptau(x)$, for all $x\in M$, and {\it $P$-central} if $y\xi=\xi y$, for all $y\in P$.
A net of vectors $\xi_n\in\mathcal K$ is called {\it $P$-almost central} if $\|y\xi_n-\xi_n y\|\rightarrow 0$, for all $y\in P$.

The analogy between representations and bimodules  led to von Neumann algebraic analogues of various representation theoretic properties of groups, including amenability and property (T):

\begin{definition} Let $(M,\uptau)$ be a tracial von Neumann algebra, and $P,Q\subset M$ be subalgebras.
\begin{enumerate}
\item We say that $M$ is {\it amenable} if there exists a net $\xi_n\in L^2(M)\bar{\otimes}L^2(M)$ of tracial, $M$-almost central vectors \cite{Po86}. We say that {\it $Q$ is amenable relative to $P$} if there exists a net $\xi_n\in L^2(M)\bar{\otimes}_PL^2(M)$ of tracial, $Q$-almost central vectors \cite{OP07}. 
\item We say that $M$ has {\it property (T)} if any Hilbert $M$-bimodule without $M$-central vectors does not admit a net of $M$-almost central unit vectors \cite{CJ83}. We  say that $P\subset M$ has the {\it relative property (T)} if any Hilbert $M$-bimodule without $P$-central vectors does not admit a net of $M$-almost central, tracial vectors \cite{Po01b}.
\end{enumerate}
\end{definition}

\subsection{Popa's intertwining-by-bimodules} In \cite{Po01b,Po03}, Popa developed a powerful technique for showing unitary conjugacy of subalgebras of a tracial von Neumann algebra.

\begin{theorem}[Popa,\cite{Po01b,Po03}]\label{corner}
If $P, Q$ are von Neumann subalgebras of a separable tracial von Neumann algebra $(M,\uptau)$, then the following are equivalent:
\begin{enumerate}
\item There is no sequence of unitaries $u_n\in P$ satisfying $\|E_Q(au_nb)\|_2\rightarrow 0$, for all $a,b\in M$.
\item There are non-zero projections $p\in P,q\in Q$, a $*$-homomorphism $\theta:pPp\rightarrow qQq$, and a non-zero partial isometry $v\in qMp$ such that $\theta(x)v=vx$, for all $x\in pPp$.

\end{enumerate} 
If these conditions hold, we say that {\bf a corner of $P$ embeds into $Q$}. 

Moreover, if $P$ and $Q$ are Cartan subalgebras of $M$ and a corner of $P$ embeds into a corner of $Q$, then there is a unitary $u\in M$ such that $P=uQu^*$.
\end{theorem}

\section{Popa's deformation/rigidity theory}\label{3}
\subsection{Deformations}
Since its introduction in the early 2000's, Popa's deformation/rigidity theory has had a transformative impact on the theory of von Neumann algebras. 
The theory builds on Popa's  innovative idea of using the deformations of a  II$_1$ factor to locate its rigid subalgebras.
Before illustrating this principle with several examples, let us make precise the notion of a deformation.

\begin{definition}
A {\it deformation} of the identity of a tracial von Neumann algebra $(M,\uptau)$ is a sequence of unital, trace preserving, completely positive maps $\phi_n:M\rightarrow M$ satisfying $$\|\phi_n(x)-x\|_2\rightarrow 0,\;\;\;\text{for all $x\in M$}.$$
\end{definition} 

A linear map $\phi:M\rightarrow M$ is called {\it completely positive} if the map $\phi^{(m)}:\mathbb M_m(M)\rightarrow\mathbb M_m(M)$ given by $\phi^{(m)}([x_{i,j}])=[\phi(x_{i,j})]$ is positive, for all $m\geq 1$. Note that any unital, trace preserving, completely positive map $\phi:M\rightarrow M$ extends to a contraction $\phi:L^2(M)\rightarrow L^2(M)$.

\begin{remark}\label{mall}
Deformations arise naturally from continuous families of automorphisms of larger von Neumann algebras.
To be precise, let $(\tilde M,\tilde\uptau)$ be a tracial von Neumann algebra containing $M$ such that  ${\tilde\uptau}_{|M}=\uptau$. Assume that $(\theta_t)_{t\in\mathbb R}$ is a pointwise $\|.\|_2$-continuous family of  trace preserving automorphisms of $\tilde M$ with $\theta_0=\text{id}$.
Then  $\phi_n:=E_{M}\circ\theta_{t_n}:M\rightarrow M$ defines a deformation of $M$, for any sequence $t_n\rightarrow 0$. Abusing notation, such pairs $(\tilde M,(\theta_t)_{t\in\mathbb R})$ are also called deformations of $M$.
\end{remark}

Next, we proceed to give three examples of deformations.
For a comprehensive list of examples, we refer the reader to \cite[Section 3]{Io12b}.

\begin{example}\label{1}
First, let $\Gamma$ be a countable group and $\varphi_n:\Gamma\rightarrow \mathbb C$ be a sequence of positive definite functions such that $\varphi_n(e)=1$, for all $n$, and $\varphi_n(g)\rightarrow 1$, for all $g\in\Gamma$. Then $\phi_n(u_g)=\varphi_n(g)u_g$ defines a deformation of the group algebra $L(\Gamma)$ and  $\phi_n(au_g)=\varphi_n(g)au_g$ defines a deformation of any group measure space algebra $L^{\infty}(X)\rtimes\Gamma$. 

If $\Gamma$ has Haagerup's property \cite{Ha79}, then there is such a sequence $\varphi_n:\Gamma\rightarrow\mathbb C$ satisfying $\varphi_n\in c_0(\Gamma)$, for all $n$.
When applied to $\Gamma=\text{SL}_2(\mathbb Z)$, the above procedure gives a deformation of the II$_1$ factor $M=L^{\infty}(\mathbb T^2)\rtimes$ SL$_2(\mathbb Z)$ which is compact relative to $L^{\infty}(\mathbb T^2)$. 
This fact was a crucial ingredient in Popa's proof that $M$ has a trivial fundamental group \cite{Po01b}. 
 
\end{example}

\begin{example}\label{2} Second, let $\Gamma\curvearrowright (X,\mu):=([0,1]^{\Gamma},{\bf Leb}^{\Gamma})$ be the Bernoulli action of a countable group $\Gamma$. In \cite{Po01a,Po03}, Popa discovered that Bernoulli actions have a remarkable deformation property, called {\it malleability}: there is a continuous family of automorphisms $(\alpha_t)_{t\in\mathbb R}$  of the product space  $X\times X$ which commute with diagonal action of $\Gamma$ and satisfy $\alpha_0=\text{id}$ and $\alpha_1(x,y)=(y,x)$.

To see this, first construct a continuous family of automorphisms  $(\alpha_t^0)_{t\in\mathbb R}$ of the probability space $[0,1]\times [0,1]$ such that $\alpha_0^0=\text{id}$ and $\alpha_1^0(x,y)=(y,x)$. For example, we can take $$\alpha_t^0(x,y)=\begin{cases} (x,y),\;\;\;\text{if $|x-y|\geq t$}\\ (y,x),\;\;\;\text{if $|x-y|<t$.}\end{cases}$$

Then identify $X\times X=([0,1]\times [0,1])^{\Gamma}$ and define $\alpha_t((x_g)_{g\in\Gamma})=(\alpha_t^0(x_g))_{g\in\Gamma}$.

The automorphisms $(\alpha_t)_{t\in\mathbb R}$ of $X\times X$ give rise to a deformation of $M:=L^{\infty}(X)\rtimes\Gamma$, as follows.
Since $\alpha_t$ commutes with the diagonal action of $\Gamma$ on $X\times X$, the formula
$\theta_t(au_g)=(a\circ\alpha_t^{-1})u_g$ defines a trace preserving automorphism of $\tilde M:=L^{\infty}(X\rtimes X)\rtimes\Gamma$. Thus, $(\theta_t)_{t\in\mathbb R}$ is a continuous family of automorphisms of $\tilde M$ such that $\theta_0=\text{id}$ and $\theta_1(a\otimes b)=b\otimes a$, for all $a,b\in L^{\infty}(X)$. 
Since $M$ embeds into $\tilde M$ via the map $au_g\mapsto (a\otimes 1)u_g$, one obtains a deformation of $M$ (see Remark \ref{mall}). 
\end{example}

\begin{example}\label{3}
Finally, we recall  from \cite{Po86,Po06c} the construction of a malleable deformation for the free group factors, $L(\mathbb F_n)$. For simplicity, we consider the case $n=2$ and put $M=L(\mathbb F_n)$. Denote by $a_1,a_2,b_1,b_2$ the generators of $\mathbb F_4$, and view $\mathbb F_2$ as the subgroup of $\mathbb F_4$ generated by $a_1,a_2$.
This gives an embedding of $M$  into $\tilde M=L(\mathbb F_4)$. If we see $b_1$ and $b_2$ as unitary elements of $\tilde M$, then we can find self-adjoint operators $h_1$ and $h_2$ such that $b_1=\exp(i h_1)$ and $b_2=\exp(i h_2)$.
One can now define a $1$-parameter group of automorphism $(\theta_t)_{t\in\mathbb R}$ of $\tilde M$ as follows:
$$\theta_t(a_1)=\exp(ith_1)a_1,\;\;\;\theta_t(a_2)=\exp(ith_2)a_2,\;\;\;\theta_t(b_1)=b_1,\;\;\;\text{and}\;\;\theta_t(b_2)=b_2.$$
\end{example}

\subsection{Deformation vs. rigidity}\label{defrig} We will now explain briefly and informally how Popa used these deformations to prove structural results for subalgebras satisfying various rigidity properties. 

A main source of rigidity is provided by the relative property (T). Indeed, assume that $P\subset M$ is an inclusion with the relative property (T).
The correspondence between completely positive maps and bimodules (see, e.g., \cite[Section 2.1]{Po06b}) implies that any deformation $\phi_n:M\rightarrow M$ must converge uniformly to the identity in $\|.\|_2$ on the unit ball of $P$ \cite{Po01b}. In particular, for any large enough $n_0\geq 1$ one has that $$(\text{i})\;\;\;\;\|\phi_{n_0}(u)-u\|_2\leq 1/2\;\;\;\text{for all unitaries $u\in P$.}$$
In \cite{Po01b,Po03}, this analytical condition is combined with the intertwining-by-bimodule technique to deduce $P$ can be unitarily conjugate into a distinguished subalgebra $Q\subset M$.

First, suppose that $M=L^{\infty}(\mathbb T^2)\rtimes\Gamma$, where $\Gamma=\text{SL}_2(\mathbb Z)$, and denote $Q=L^{\infty}(\mathbb T^2)$.
Since $\Gamma$ has Haagerup's property, we can find positive definite functions $\varphi_n\in c_0(\Gamma)$ such that $\varphi_n\rightarrow 1$ pointwise.  As in Example \ref{1}, we obtain a deformation $\phi_n:M\rightarrow M$ given by $\phi_n(au_g)=\varphi_n(g)au_g$. Consider the Fourier decomposition $u=\sum_{g\in\Gamma}E_Q(uu_g^*)u_g$ of a unitary $u$ belonging to the subalgebra $P\subset M$ with the relative property (T). The specific formula of $\phi_{n_0}$ allows one to rewrite $(\text{i})$ as $$(\text{ii})\;\;\;\;\sum_{g\in\Gamma}|\varphi_{n_0}(g)-1|^2\;\|E_Q(uu_g^*)\|_2^2\leq 1/4\;\;\;\text{for all unitaries $u\in P$.}$$
Since $\varphi_{n_0}\in c_0(\Gamma)$, we have that $|\varphi_{n_0}(g)-1|^2\geq 1/2$, for all $g\in\Gamma$ outside a finite subset $F$.
Taking into account that $\sum_{g\in\Gamma}\|E_Q(uu_g^*)\|_2^2=\uptau(u^*u)=1$ and using $(\text{ii})$, one concludes that $$(\text{iii})\;\;\;\;\sum_{g\in F}\|E_Q(uu_g^*)\|_2^2\geq 1/2\;\;\;\text{for all unitaries $u\in P$}.$$
As a corollary, $P$ does not admit a sequence a unitaries satisfying condition (1) of Proposition \ref{corner}. In other words, a corner of $P$ embeds into $Q$. If $P$ is a Cartan subalgebra of $M$, then the moreover part of Proposition \ref{corner} implies that $P$ must be unitarily conjugate to $Q$ \cite{Po01b}.

 Second, suppose that $M=L^{\infty}(X)\rtimes\Gamma$, where $\Gamma\curvearrowright (X,\mu)=([0,1]^{\Gamma},{\bf Leb}^{\Gamma})$ is the Bernoulli action. Denote $Q=L(\Gamma)$ and let $(\tilde M,(\theta_t)_{t\in\mathbb R})$ the deformation introduced in Example \ref{2}.
Then the deformation $\phi_n:=E_M\circ\theta_{1/2^n}$ converges uniformly to the identity in $\|.\|_2$ on the unit ball of any subalgebra $P\subset M$ with the property (T). It is immediate that the same must be true for $\theta_{1/2^n}$. In particular, for any large enough $n_0\geq 1$ one has that $$(\text{iv})\;\;\;\;\|\theta_{1/2^{n_0}}(u)-u\|_2\leq 1/2\;\;\;\text{for all unitaries $u\in P$.}$$

This implies the existence of a non-zero element $v\in\tilde M$ satisfying $\theta_{1/2^{n_0}}(u)v=vu$ for all unitaries $u\in P$. By employing a certain symmetry of the deformation $(\theta_t)_{t\in\mathbb R}$, Popa proved that the same holds for $\theta_1$.  Using the formula of the restriction $\theta_1$ to $M$, it follows that a corner of $P$ embeds into $Q$. Assuming that $\Gamma$ is icc,  Popa concludes that $P$ can be unitarily conjugate into $Q$ \cite{Po03}. 

The rigidity considered in \cite{Po01b,Po03} is due to (relative) property (T) assumptions.  In  \cite{Po06a,Po06c}, Popa discovered a less restrictive from of rigidity, arising from the presence of subalgebras with non-amenable relative commutant. To illustrate this, assume the context from Example \ref{3}: $M=L(\mathbb F_2)\subset\tilde M=L(\mathbb F_4)$ and $(\theta_t)_{t\in\mathbb R}$ is the $1$-parameter group of automorphisms of $\tilde M$ defined therein. The deformation $(\tilde M, (\theta_t)_{t\in\mathbb R})$ has two crucial properties:
\begin{enumerate}[label=(\alph*)]
\item the $M$-bimodule $L^2(\tilde M)\ominus L^2(M)$ is isomorphic to a multiple of $L^2(M)\bar{\otimes}L^2(M)$.
\vskip 0.05in
\item the contraction $E_M\circ\theta_t:L^2(M)\rightarrow L^2(M)$ is a compact operator, for any $t>0$.
\end{enumerate}

Recovering Ozawa's solidity theorem \cite{Oz03} in the case of the free group factors, Popa proved in \cite{Po06c} that the relative commutant $P'\cap M$ is amenable, for any diffuse subalgebra $P\subset M$. Let us explain how  (a) and (b) are combined in \cite{Po06c} to deduce the following weaker statement: there is no diffuse subalgebra $P\subset M$ which commutes with a non-amenable II$_1$ subfactor $N\subset M$. 

Assume that we can find such  commuting subalgebras $P,N\subset M$. The non-amenability of $N$ leads to a spectral gap condition for its coarse bimodule \cite{Co75b}: there is a finite set $F\subset N$ such that $$\text{(v)}\;\;\;\;\|\xi\|\leq\sum_{x\in F}\|x\cdot\xi-\xi\cdot x\|\;\;\;\text{for all vectors $\xi\in L^2(N)\bar{\otimes}L^2(N)$}.$$ 

The spectral gap condition is then used  to establish the following rigidity property for $P$: the deformation $E_{M}\circ\theta_t:M\rightarrow M$ converges uniformly on the unit ball of $P$.
One first notes that  (a) implies that (v) holds for every $\xi\in\tilde M$ with $E_M(\xi)=0$. Thus, we may take $\xi=\theta_t(u)-E_M(\theta_t(u))$, for any unitary $u\in P$ and $t>0$. Using that $u$ commutes with every $x\in F$, one derives that
\begin{align*}
\text{(vi)}\;\;\;\;\|x\cdot\xi-\xi\cdot x\|&\leq \|x\theta_t(u)-\theta_t(u)x\|_2\\ &=\|\theta_{-t}(x)u-u\theta_{-t}(x)\|_2\\&=\|(\theta_{-t}(x)-x)u-u(\theta_{-t}(x)-x)\|_2\\&\leq 2\|\theta_{-t}(x)-x\|_2.
\end{align*}
Combining (v) and (vi) gives that if $t>0$ is chosen small enough then $\|\theta_t(u)-E_{M}(\theta_t(u))\|_2\leq 1/2$, for all unitaries $u\in P$. Since $P$ is assumed diffuse it contains a sequence of unitaries $u_n$ converging weakly to $0$. The compactness of $E_M\circ\theta_t$  gives that $\|E_{M}(\theta_t(u_n))\|_2\rightarrow 0$, leading to a contradiction.

\section{Uniqueness of Cartan subalgebras and W$^*$-superrigidity of Bernoulli actions}

\subsection{Uniqueness of Cartan subalgebras} 
 The first result showing uniqueness, up to unitary conjugacy, of arbitrary Cartan subalgebras, was proved by Ozawa and Popa: 

\begin{theorem}[Ozawa, Popa, \cite{OP07}]\label{OP07} Let $\mathbb F_n\curvearrowright (X,\mu)$ be a free ergodic profinite p.m.p. action of $\mathbb F_n$, for $n\geq 2$.
Then $M:=L^{\infty}(X)\rtimes\mathbb F_n$ has a unique Cartan subalgebra, up to unitary conjugacy: if $P\subset M$ is any Cartan subalgebra, then $P=uL^{\infty}(X)u^*$, for some unitary element $u\in M$.
\end{theorem}

The proof of Theorem \ref{OP07} relies on the {\it complete metric approximation property} (CMAP) of $\mathbb F_n$, used as a weak form of a deformation of the II$_1$ factor $L(\mathbb F_n)$. Recall that a countable group $\Gamma$ has the CMAP \cite{Ha79} if there exists a sequence of finitely supported functions $\varphi_k:\Gamma\rightarrow\mathbb C$ such that $\varphi_k(g)\rightarrow 1$, for all $g\in\Gamma$, and the linear maps $\phi_k:L(\Gamma)\rightarrow L(\Gamma)$ given by $\phi_k(u_g)=\varphi_k(g)u_g$, for all $g\in\Gamma$, satisfy $\limsup_k\|\phi_k\|_{\text{cb}}=1$.
If the last condition is weakened by assuming instead that $\limsup_k\|\phi_k\|_{\text{cb}}<\infty$, then $\Gamma$ is called {\it weakly amenable} \cite{CH88}. 
Here, $\|\phi_k\|_{\text{cb}}$ denotes the completely bounded norm of $\phi_k$.

Since $\mathbb F_n$ has the CMAP \cite{Ha79} and the action $\mathbb F_n\curvearrowright X$ is profinite, the II$_1$ factor $M$ also has the CMAP: there exists a sequence of finite rank completely bounded maps $\phi_k:M\rightarrow M$ such that $\|\phi_k(x)-x\|_2\rightarrow 0$, for all $x\in M$, and $\limsup_k\|\phi_k\|_{\text{cb}}=1$. 

Let $P\subset M$ be an arbitrary diffuse amenable subalgebra and denote by $\mathcal G$ its normalizer in $M$. If the conjugation action $\mathcal G\curvearrowright P$ happens to be compact, i.e. the closure of $\mathcal G$ inside Aut$(P)$ is compact, then the unitary representation $\mathcal G\curvearrowright L^2(P)$ is a direct sum of finite dimensional representations.  This provides many vectors $\xi\in L^2(P)\bar{\otimes}L^2(P)$ that are invariant under the diagonal action of $\mathcal G$.  Indeed, if $\eta_1,...,\eta_d$ is an orthonormal basis of any finite dimensional $\mathcal G$-invariant subspace of $L^2(P)$, then $\xi=\sum_{i=1}^d\eta_i\otimes\eta_i^*$ has this property.

 Ozawa and Popa made the fundamental discovery that since $M$ has the CMAP, the action $\mathcal G\curvearrowright P$ is {\it weakly compact} (although it is typically not compact). More precisely, they showed the existence of a net of vectors $\xi_k\in L^2(P)\bar{\otimes}L^2(P)$ which are almost invariant under the diagonal action of $\mathcal G$.  In the second part of the proof, they combined the weak compactness property, with a malleable deformation of $M$ analogous to the one from Example \ref{3} and a spectral gap rigidity argument. Thus, they concluded that either a corner of $P$ embeds into $L^{\infty}(X)$ or the von Neumann algebra generated by $\mathcal G$ is amenable. 
If $P\subset M$ is a Cartan subalgebra, then the first condition must hold since $M$ is not amenable. By Theorem \ref{corner}, this forces that $P$ is unitarily conjugate to $L^{\infty}(X)$.

Theorem \ref{OP07} is restricted both by the class of groups and the family of actions it applies to. 
Ozawa \cite{Oz10} showed that the weak compactness property established in the proof of Theorem \ref{OP07} more generally holds for profinite actions of weakly amenable groups. 
Building on this and the weak amenability of hyperbolic groups \cite{Oz07}, Chifan and Sinclair extended Theorem \ref{OP07} to all non-elementary hyperbolic groups $\Gamma$ \cite{CS11}.  A conceptual novelty of their approach was the usage of quasi-cocyles (rather than cocycles \cite{PS09,Si10}) to build 
deformations. 

Soon after,  Popa and Vaes obtained uniqueness results of unprecedented generality in \cite{PV11,PV12}. Generalizing 
\cite{OP07,CS11}, they showed that Theorem \ref{OP07} holds for arbitrary actions of free groups and hyperbolic groups:

\begin{theorem}[Popa, Vaes, \cite{PV11,PV12}]\label{PV} Let $\Gamma\curvearrowright (X,\mu)$ be a free ergodic p.m.p. action of a countable group $\Gamma$. Assume either that $\Gamma$ is weakly amenable and admits an unbounded cocycle into a non-amenable mixing orthogonal representation, or $\Gamma$ is non-elementary hyperbolic.	
Then $L^{\infty}(X)\rtimes\Gamma$ has a unique Cartan subalgebra, up to unitary conjugacy.
\end{theorem} 

This result covers any weakly amenable group $\Gamma$
with a positive first $\ell^2$-Betti number, $\beta_1^{(2)}(\Gamma)>0$. Indeed, the latter holds if and only if $\Gamma$ is non-amenable and admits an unbounded cocycle into its left regular representation \cite{BV97,PT07}. 

Theorem \ref{PV} led to the resolution of the group measure space analogue of the famous, still unsolved, {\it free group factor problem} which asks whether $L(\mathbb F_n)$ and $L(\mathbb F_m)$ are isomorphic or not, for $n\not=m$.
More precisely, Popa and Vaes showed in \cite{PV11} that if $2\leq n,m\leq\infty$ and $n\not=m$, then for any free ergodic p.m.p. actions $\mathbb F_n\curvearrowright (X,\mu)$ and $\mathbb F_m\curvearrowright (Y,\nu)$ one has: $$L^{\infty}(X)\rtimes\mathbb F_n\not\cong L^{\infty}(Y)\rtimes\mathbb F_m.$$
If these factors were isomorphic, then Theorem \ref{PV} would imply that the actions $\mathbb F_n\curvearrowright X$, $\mathbb F_m\curvearrowright Y$ are orbit equivalent. However, it was shown by Gaboriau \cite{Ga99,Ga01} that free groups of different ranks do not admit orbit equivalent free actions.

In the setting of Theorem \ref{PV}, let $P\subset L^{\infty}(X)\rtimes\Gamma$ be a diffuse amenable subalgebra and denote by $\mathcal G$ its normalizer.
The weak amenability of $\Gamma$ is used to show that, roughly speaking, the action $\mathcal G\curvearrowright P$ is weakly compact relative to $L^{\infty}(X)$. When $P$ is a Cartan subalgebra, in  combination with the deformations obtained from cocycles or quasi-cocycles of $\Gamma$, one concludes that $P$ is unitarily conjugate to $L^{\infty}(X)$. 

Now, assume that $\Gamma$ is a non-abelian free group or, more generally,  a non-elementary hyperbolic group. In this context, Popa and Vaes proved a deep, important dichotomy for arbitrary diffuse amenable subalgebras $P$ of arbitrary tracial crossed products $B\rtimes\Gamma$: either a corner of $P$ embeds into $B$, or otherwise the von Neumann algebra generated by the normalizer of $P$ is amenable relative to $B$ \cite{PV11,PV12}.    
This property of such groups $\Gamma$, called {\it relative strong solidity}, has since found a number of impressive applications. In particular, it was used in \cite{Io12a} to provide the first class of non-weakly amenable groups satisfying Theorem \ref{PV}:

\begin{theorem}\emph{\cite{Io12a}}\label{Io12}
Let $\Gamma=\Gamma_1*_{\Sigma}\Gamma_2$ be an amalgamated free product group. Assume that $[\Gamma_1:\Sigma]\geq 2$, $[\Gamma_2:\Sigma]\geq 3$, and $\cap_{i=1}^n g_i\Sigma g_i^{-1}=\{e\}$, for some elements $g_1,...,g_n\in\Gamma$. Let $\Gamma\curvearrowright (X,\mu)$ be a free ergodic p.m.p. action.
Then $L^{\infty}(X)\rtimes\Gamma$ has a unique Cartan subalgebra, up to unitary conjugacy.
\end{theorem}

Theorem \ref{Io12} in particular applies to any free product group $\Gamma=\Gamma_1*\Gamma_2$ with $|\Gamma_1|\geq 2$, $|\Gamma_2|\geq 3$.
Since such groups have a positive first $\ell^2$-Betti number, Theorems \ref{PV} and \ref{Io12} both provide supporting evidence for the following general conjecture:

{\bf Problem I}. {\it Let $\Gamma$ be a countable group with $\beta_1^{(2)}(\Gamma)>0$. Prove that $L^{\infty}(X)\rtimes\Gamma$ has a unique Cartan subalgebra, up to unitary conjugacy, for any free ergodic p.m.p. action $\Gamma\curvearrowright (X,\mu)$.}

A weaker version of Problem I asks to prove that $L^{\infty}(X)\rtimes\Gamma$ has a unique group measure space Cartan subalgebra. This has been confirmed by Chifan and Peterson in \cite{CP10} under the additional assumption that $\Gamma$ has a non-amenable subgroup with the relative property (T) (see also \cite{Va10b}). A positive answer was also obtained  in \cite{Io11a,Io11b} when the action $\Gamma\curvearrowright (X,\mu)$ is rigid or profinite.

If Problem I or its weaker version were solved, then as $\ell^2$-Betti numbers of groups are invariant under orbit
equivalence \cite{Ga01}, it would follow that $\beta_1^{(2)}(\Gamma)$  is an isomorphism invariant of $L^{\infty}(X)\rtimes\Gamma$. This would provide a computable invariant for the class of group measure space II$_1$ factors. 

\subsection{Non-uniqueness of Cartan subalgebras} 
If $\Gamma$ is any group as in Theorems \ref{PV} and \ref{Io12}, then $L^{\infty}(X)\rtimes\Gamma$ has a unique Cartan subalgebras, for any free ergodic p.m.p. action of $\Gamma$, while $L(\Gamma)$  and $L(\Gamma)\bar{\otimes}N$ do not have Cartan subalgebras, for any II$_1$ factor $N$. This provides several large families of II$_1$ factors with at most one Cartan subalgebra. However, in the non-uniqueness regime, little is known about the possible cardinality of the set of Cartan subalgebras of a II$_1$ factor.

Although the hyperfinite II$_1$ factor $R$ admits uncountably many Cartan subalgebras up to unitary conjugacy \cite{Pa85}, any two Cartan subalgebras are
conjugated by an automorphism of $R$ \cite{CFW81}. The first class of examples of II$_1$ factors admitting two
Cartan subalgebras that are not conjugated by an automorphism was given by
Connes and Jones in \cite{CJ81}. A second class of examples of such II$_1$ factors where the two Cartan
subalgebras are explicit was found by Ozawa and Popa \cite{OP08} (see also
\cite{PV09}). A class of II$_1$ factors $M$ whose Cartan subalgebras cannot be concretely classified up to unitary conjugacy or up to conjugation by an automorphism of $M$ was then introduced in \cite{SV11}. Most recently, a family of II$_1$ factors whose all group measure space Cartan subalgebras can be described explicitly was constructed in \cite{KV15}. In particular, the authors give examples of II$_1$ factors having exactly two group measure space Cartan subalgebras, up to unitary conjugacy, and a prescribed number of group measure space Cartan subalgebras, up to conjugacy with an automorphism. However, there are currently no known examples of II$_1$ having precisely $n\geq 2$ arbitrary Cartan subalgebras:

\vskip 0.05in
{\bf Problem II}. {\it Given an integer $n\geq 2$, find II$_1$ factors $M$ which have exactly $n$ Cartan subalgebras, up to unitary conjugacy (or up to conjugacy with an automorphism).}

\subsection{W$^*$-superrigidity of Bernoulli actions} Popa's strong rigidity theorem \cite{Po04a} led to the natural conjecture that Bernoulli actions of icc property (T)  groups are W$^*$-superrigid. In this section we discuss the solution of this conjecture.

\begin{theorem}[Ioana, \cite{Io10}]\label{Io10} Let $\Gamma$ be an icc property (T) group and $(X_0,\mu_0)$ be a non-trivial probability space. 
Then the Bernoulli action $\Gamma\curvearrowright (X,\mu)=(X_0^{\Gamma},\mu_0^{\Gamma})$ is W$^*$-superrigid.
\end{theorem}

The proof of Theorem \ref{Io10} relies on a general strategy for analyzing group measure space decompositions of II$_1$ factors \cite{Io10}. 
Let $M=L^{\infty}(X)\rtimes\Gamma$ be the II$_1$ factor arising from a ``known" free ergodic p.m.p. action $\Gamma\curvearrowright (X,\mu)$.
Assume that we can also decompose $M=L^{\infty}(Y)\rtimes\Lambda$, for some ``mysterious" free ergodic p.m.p. action $\Lambda\curvearrowright (Y,\nu)$.

The new group measure space decomposition of $M$ gives rise to an embedding $\Delta:M\rightarrow M\bar{\otimes}M$ defined by $\Delta(bu_h)=u_h\otimes u_h$, for all $b\in L^{\infty}(X)$ and $h\in\Lambda$.
This embedding has been introduced in \cite{PV09} were it was used to transfer rigidity properties through W$^*$-equivalence.
The strategy of \cite{Io10} is to first prove a classification of all embeddings of $M$ into $M\bar{\otimes}M$, and then apply it to $\Delta$. This provides a relationship between the actions $\Gamma\curvearrowright X,\Lambda\curvearrowright Y$ which is typically stronger than the original W$^*$-equivalence and which, ideally, can be exploited to show that the actions are orbit equivalent or even conjugate.

In the case $\Gamma\curvearrowright (X,\mu)$ is a Bernoulli action of a property (T) group, a classification of all possible embeddings $\theta:M\rightarrow M\bar{\otimes}M$ was obtained in \cite{Io10}. This classification is precise enough so that when combined with the above strategy it implies that $L^{\infty}(X)$ is the unique group measure space Cartan subalgebra of $M$, up to unitary conjugacy. Since the action $\Gamma\curvearrowright (X,\mu)$ is OE-superrigid by a theorem of Popa \cite{Po05}, it follows that the action is also W$^*$-superrigid. 

The main novelty of \cite{Io10} is a structural result for abelian subalgebras $D$ on $M$ that are normalized by a sequence of unitary elements $u_n\in L(\Gamma)$ converging weakly to $0$.  Under this assumption, it is shown that $D$ and its relative commutant $D'\cap M$ can be essentially unitarily conjugated into either $L(\Gamma)$ or $L^{\infty}(X)$.  An analogous dichotomy is proven in \cite{Io10} for abelian subalgebras $D\subset M\bar{\otimes}M$ which are normalized by ``many" unitary elements from $L(\Gamma)\bar{\otimes}L(\Gamma)$.
This is applied to study embeddings $\theta:M\rightarrow M\bar{\otimes}M$ as follows. Since $\theta(L(\Gamma))$ is a property (T) subalgebra of $M\bar{\otimes}M$, by adapting the arguments described in Section \ref{defrig},  we may assume that $\theta(L(\Gamma))\subset L(\Gamma)\bar{\otimes}L(\Gamma)$. As a consequence, $D=\theta(L^{\infty}(X))$ is normalized by the group $\theta(\Gamma)\subset L(\Gamma)\bar{\otimes}L(\Gamma)$, and thus the dichotomy can be applied to $D$.  

By Theorem \ref{Io10}, II$_1$ factors arising from Bernoulli actions of icc property (T) groups have a unique group measure space Cartan subalgebra. The following problem proposed by Popa asks to prove that the same holds for arbitrary non-amenable groups and general Cartan subalgebras: 

\vskip 0.05in
{\bf Problem III.} {\it Let $\Gamma\curvearrowright (X,\mu)=(X_0^{\Gamma},\mu_0^{\Gamma})$ be a Bernoulli action of a non-amenable group $\Gamma$. 
Then $L^{\infty}(X)\rtimes\Gamma$ has a unique Cartan subalgebra, up to unitary conjugacy.}

A positive answer to this problem would imply that if two Bernoulli actions of non-amenable groups are W$^*$-equivalent, then they are orbit equivalent. By  \cite{Po05,Po06b},  for Bernoulli actions of groups in a large class (containing all infinite property (T) groups), orbit equivalence implies conjugacy.  However, this does not hold for arbitrary non-amenable groups. Indeed, if $n\geq 2$, then the Bernoulli actions $\mathbb F_n\curvearrowright (X_0^{\mathbb F_n},\mu_0^{\mathbb F_n})$ of $\mathbb F_n$ are completely classified up to conjugacy by the entropy of base space $(X_0,\mu_0)$ \cite{Bo08}; on the other hand, all Bernoulli actions of $\mathbb F_n$ are orbit equivalent \cite{Bo09}.
 
\section{Orbit equivalence rigidity}

Pioneering OE rigidity results were obtained by Zimmer for actions of higher rank semisimple Lie groups and their lattices  by using his influential cocycle superrigidity theorem \cite{Zi80,Zi84}.  
Deducing OE rigidity results from cocycle superrigidity theorems has since become a paradigm in the area. 
An illustration of this is Furman's proof that ``generic" free ergodic p.m.p. actions of higher rank lattices, including the actions   SL$_n(\mathbb Z)\curvearrowright(\mathbb T^n,\text{\bf Leb})$ for $n\geq 3$, are virtually  OE-superrigid. Since then, numerous striking OE superrigidity results have been discovered in \cite{Po05,Po06a,Ki06,Io08,PV08b,Fu09,Ki09,PS09,TD14,Io14,CK15,Dr15, GITD16}.

\subsection{Cocycle superrigidity} Many of these results have been obtained by applying techniques and ideas from Popa's deformation/rigidity theory. In all of these cases, one proves that much more than being OE superrigid, the actions in question are cocycle superrigid  \cite{Po05,Po06a,Io08,PV08b,Fu09,PS09,TD14,Io14,Dr15, GITD16}.
These developments were trigerred by Popa's discovery of a striking new cocycle superrigidity phenomenon:

\begin{theorem}[Popa, \cite{Po05,Po06a}]\label{Po05}
Assume that $\Gamma$ is either an infinite property (T) group or the product $\Gamma_1\times\Gamma_2$ of an infinite group and a non-amenable group.
Let $\Gamma\curvearrowright (X,\mu):=(X_0^{\Gamma},\mu_0^{\Gamma})$ be the Bernoulli action, where $(X_0,\mu_0)$ is a non-trivial standard probability space. Let $\Lambda$ be a countable group. 
Then any cocycle $w:\Gamma\times X_0^{\Gamma}\rightarrow\Lambda$ is cohomologous to a group homomorphism $\delta:\Gamma\rightarrow\Lambda$.
\end{theorem}

Theorem \ref{Po05} more generally applies to cocycles with values into $\mathcal U_{\text{fin}}$ groups $\Lambda$, i.e.,  isomorphic copies of closed subgroups of the unitary group of a separable II$_1$ factor.

Recall that if $\Lambda$ is a Polish group, then a measurable map $w:\Gamma\times X\rightarrow\Lambda$ is called a {\it cocycle} if  it satisfies the identity $w(g_1g_2,x)=w(g_1,g_2\cdot x)w(g_2,x)$ for all $g_1,g_2\in\Gamma$ and almost every $x\in X$. Two cocycles $w_1,w_2:\Gamma\times X\rightarrow\Lambda$ are called {\it cohomologous} if there exists a measurable map $\varphi:X\rightarrow\Lambda$ such that $w_2(g,x)=\varphi(g\cdot x)w_1(g,x)\varphi(x)^{-1}$, for all $g\in\Gamma$ and almost every $x\in X$. Any group homomorphism $\delta:\Gamma\rightarrow\Lambda$ gives rise to a constant cocycle, $w(g,x):=\delta(g)$. 

\begin{remark} If $\alpha:(X,\mu)\rightarrow (Y,\nu)$ is an OE between $\Gamma\curvearrowright (X,\mu)$ and a free p.m.p. action $\Lambda\curvearrowright (Y,\nu)$, then the map $w:\Gamma\times X\rightarrow \Lambda$ uniquely determined by the formula $\alpha(gx)=w(g,x)\cdot\alpha(x)$ is a cocycle, called the {\it Zimmer cocycle}.
Assume that $w$ is cohomologous to a homomorphism, i.e. $w(g,x)=\varphi(g\cdot x)\delta(g)\varphi(x)^{-1}$. Then the map $\tilde\alpha:X\rightarrow Y$ given by $\tilde\alpha(x):=\varphi(x)^{-1}\cdot\alpha(x)$ satisfies $\tilde\alpha(g\cdot x)=\delta(g)\cdot\tilde\alpha(x)$. This can be often used to conclude that the actions are conjugate, e.g., if $\Gamma$ is icc and its action is free and weakly mixing \cite{Po05}. Consequently, the actions from Theorem \ref{Po05} are OE-superrigid whenever $\Gamma$ is icc.
\end{remark}

The proof of Theorem \ref{Po05} relies on the malleability of Bernoulli actions. Assume for simplicity that $(X_0,\mu_0)=([0,1],\text{\bf Leb})$. Let $(\alpha_t)_{t\in\mathbb R}$ be a continuous family of automorphisms of $X\times X$ commuting with the diagonal action of $\Gamma$ and satisfying $\alpha_0=\text{id}$ and $\alpha_1(x,y)=(y,x)$, as in Example \ref{2}. Then the formula $w_t(g,(x,y))=w(g,\alpha_t(x,y))$ defines a  family $(w_t)_{t\in\mathbb R}$ of cocycles for the product action $\Gamma\curvearrowright X\times X$. Note that $w_0(g,(x,y))=w(g,x)$, while $w_1(g,(x,y))=w(g,y)$.

In \cite{Po05}, Popa uses property (T) to deduce that $w_t$ is cohomologous to $w_0$, for $t>0$ small enough. The same conclusion is derived in \cite{Po06a} via a spectral gap rigidity argument (see Section \ref{defrig}). In both papers, this is then shown to imply that $w$ is cohomologous to $w_1$. Finally, using the weak mixing property of Bernoulli actions,  it is concluded that  $w$ is cohomologous to a homomorphism.
For a presentation of  the proof of Theorem \ref{Po05} in the property (T) case, see also \cite{Fu06,Va06a}.

In addition to property (T) and product groups, Theorem \ref{Po05} has been shown to hold for groups $\Gamma$ such that $L(\Gamma)$ is L$^2$-rigid in the sense of \cite{Pe06} (see \cite{PS09}) and for inner amenable non-amenable groups $\Gamma$ (see \cite{TD14}). The following problem due to Popa  remains however open:

{\bf Problem IV.} {\it Characterize the class of groups $\Gamma$ whose Bernoulli actions are cocycle superrigid (for arbitrary countable or $\mathcal U_{\text{fin}}$ ``target" groups $\Lambda$), in the sense of Theorem \ref{Po05}.}

This class is conjecturally characterised by the vanishing of the first L$^2$-Betti number. 
 Indeed, all groups $\Gamma$ known to belong to the class satisfy $\beta_1^{(2)}(\Gamma)=0$. 
On the other hand, if $\beta_1^{(2)}(\Gamma)>0$, then the Bernoulli actions of $\Gamma$ are not cocycle superrigid with $\Lambda=\mathbb T$ as the target group \cite{PS09}. 

 In \cite{Io08}, the author established a cocycle superrigidity theorem for translation actions  $\Gamma\curvearrowright (G,{\bf m}_G)$ of property (T) groups $\Gamma$ on their profinite completions $G$, see part (1) of Theorem \ref{profinite}. 
Note that these actions are in some sense the farthest from being weakly mixing or Bernoulli: the unitary representation $\Gamma\curvearrowright L^2(G)$ is a  sum of finite dimensional representations, by the Peter-Weyl theorem.

Motivated by the analogy with Theorem \ref{Po05}, it was asked in \cite{Io08} whether a version of the cocycle superrigidity theorem obtained therein holds for product groups, such as $\Gamma=\mathbb F_2\times\mathbb F_2$. The interest in this question was especially high at the time, since a positive answer combined with the work \cite{OP07} would have lead to the (then) first examples of virtually W$^*$-superrigid actions. This 	question was recently settled in \cite{GITD16}, see part (2) of Theorem \ref{profinite}. 

\begin{theorem}\label{profinite} Let $\Gamma$ and $\Delta$ be countable dense subgroups of a compact profinite group $G$.
Consider the left translation action $\Gamma\curvearrowright (G,{\bf m}_G)$ and the left-right translation action $\Gamma\times\Delta\curvearrowright (G,{\bf m}_G)$. \\
  Let $\Lambda$ be a countable group. Let $w:\Gamma\times G\rightarrow\Lambda$ and $v:(\Gamma\times\Delta)\times G\rightarrow\Lambda$ be any cocycles.
\begin{enumerate}
\item \emph{\cite{Io08}} Assume that $\Gamma$ has property (T).  Then we can find an open subgroup $G_0<G$ such that the restriction of $w$ to $(\Gamma\cap G_0)\times G_0$ is cohomologous to a homomorphism $\delta:\Gamma\cap G_0\rightarrow\Lambda$.
\vskip 0.05in

\item \emph{\cite{GITD16}} Assume that $\Gamma\curvearrowright (G,{\bf m}_G)$ has spectral gap, and $\Gamma$, $\Lambda$ are finitely generated.
Then we can find an open subgroup $G_0<G$ such that the restriction of $v$ to  \\ $[(\Gamma\cap G_0)\times(\Delta\cap G_0)]\times G_0$ is cohomologous to a homomorphism $\delta:(\Gamma\cap G_0)\times(\Delta\cap G_0)\rightarrow\Lambda$.
\end{enumerate}
\end{theorem}

Recall that the left-right translation action $\Gamma\times\Delta\curvearrowright (G,{\bf m}_G)$ is given by $(g,h)\cdot x=gxh^{-1}$.
A p.m.p. action $\Gamma\curvearrowright (X,\mu)$ is said to have {\it spectral gap} if the  representation $\Gamma\curvearrowright L^2(X)\ominus\mathbb C{\bf 1}$ does not have almost invariant vectors. A well-known result of Selberg implies that the left translation action SL$_2(\mathbb Z)\curvearrowright$ SL$_2(\mathbb Z_p)$ has spectral gap, for every prime $p$. This was recently generalized in \cite{BV10} to arbitrary non-amenable subgroups $\Gamma<$ SL$_2(\mathbb Z)$: the left translation of $\Gamma$ onto its closure $\bar{\Gamma}<$ SL$_2(\mathbb Z_p)$ has spectral gap.

Note that in \cite{Fu09}, Furman provided an alternative approach to part (1) of Theorem \ref{profinite} which applies to the wider class of compact actions. 

The proof of Theorem \ref{profinite} relies on the following criterion for untwisting cocycles $w:\Gamma\times G\rightarrow\Lambda$.
First, endow the space of such cocycles with the ``uniform" metric
$$d(w,w'):=\sup_{g\in\Gamma}\;{\bf m}_G(\{x\in G\;|\;w(g,x)\not=w'(g,x)\}).$$
Second,  since the left and right translation actions of $G$ on itself commute, $w_t(g,x)=w(g,xt)$ defines a family of cocycles $(w_t)_{t\in G}$. It is  then shown in \cite{Io08,Fu09} that if $w$ verifies the uniformity condition $d(w_t,w)\rightarrow 0$, as $t\rightarrow 1_G$, then $w$ untwists, in the sense of part (1) of Theorem \ref{profinite}.

\subsection{OE rigidity for actions of non-rigid groups}
In this section, we present a rigidity result for translation actions with spectral gap which describes precisely when two such actions are OE.

Moreover, the result also applies to the notion of Borel reducibility from descriptive set theory (see e.g. the survey \cite{Th06}). 
If $\mathcal R,\mathcal S$ are equivalence relations on standard Borel spaces $X,Y$,  we say that $\mathcal R$ is {\it Borel reducible} to $\mathcal S$ whenever there exists a Borel map $\alpha:X\rightarrow Y$ such that $(x,y)\in\mathcal R\Leftrightarrow (\alpha(x),\alpha(y))\in\mathcal S$. 
This encodes  that the classification problem associated to $\mathcal R$ is no more complicated than the classification problem associated to $\mathcal S$.

\begin{theorem}\emph{\cite{Io13}}\label{Io13}
Let $\Gamma$ and $\Lambda$ be countable dense subgroups of profinite groups $G$ and $H$.
Assume that the left translation action $\Gamma\curvearrowright (G,{\bf m}_G)$ has spectral gap.
\begin{enumerate}
\item $\Gamma\curvearrowright (G,{\bf m}_G)$ is OE to $\Lambda\curvearrowright (H,{\bf m}_H)$ iff there exist open subgroup $G_0<G,H_0<H$ and a topological isomorphism $\delta:G_0\rightarrow H_0$ such that $\delta(\Gamma\cap G_0)=\Lambda\cap H_0$ and $[G:G_0]=[H:H_0]$.
\item $\mathcal R_{\Gamma\curvearrowright G}$ is Borel reducible to $\mathcal R_{\Lambda\curvearrowright H}$ iff we can find an open subgroup $G_0<G$, a closed subgroup $H_0<H$ and a topological isomorphism $\delta:G_0\rightarrow H_0$ such that $\delta(\Gamma\cap G_0)=\Lambda\cap H_0$.
\end{enumerate}
\end{theorem}

The main novelty of this theorem lies in that there are no assumptions on the groups, but instead, all the assumptions are imposed on their actions. Thus,  the theorem applies to many natural families of actions of SL$_2(\mathbb Z)$ and the free groups $\mathbb F_n$, leading to the following:

\begin{corollary}\emph{\cite{Io13}} \label{thomas}
If $S$ and $T$ are distinct non-empty sets of primes, then the actions $$\text{SL}_2(\mathbb Z)\curvearrowright\prod_{p\in S}\text{SL}_2(\mathbb Z_p)\;\;\;\;\text{and}\;\;\;\;\text{SL}_2(\mathbb Z)\curvearrowright\prod_{p\in T}\text{SL}_2(\mathbb Z_p)$$ are not orbit equivalent, and their equivalence relations are not Borel reducible one to another.
\end{corollary}

This result settles a conjecture of Thomas (see \cite[Conj. 5.7]{Th01} and \cite[Conj. 2.14]{Th06}).
In particular, it provides a continuum of treeable countable Borel equivalence relations that are pairwise incomparable with respect to Borel reducibility. 
Note that the existence of uncountably many such equivalence relations has been first established by Hjorth in \cite{Hj08}. However, prior to Corollary \ref{thomas}, not a single example of a pair of treeable countable Borel equivalence relations such that neither is Borel reducible to the other was known.

Theorem \ref{Io13} also leads to natural concrete uncountable families of pairwise non-OE free ergodic p.m.p. actions of the non-abelian free groups, $\mathbb F_n$ (the existence of such families was proved in \cite{GP03}, while the first explicit families were found in \cite{Io06b}). 
In contrast, a non-rigidity result of Bowen shows that  any two Bernoulli actions of $\mathbb F_n$ are OE \cite{Bo09}. 

Theorem \ref{Io13} admits a version which applies to connected (rather than profinite) compact groups. More precisely, let $\Gamma<G$ and $\Lambda<H$ be countable dense subgroups of compact connected Lie groups with trivial centers. Assuming that $\Gamma\curvearrowright (G,{\bf m}_G)$ has spectral gap, it is shown in \cite{Io13} that the actions $\Gamma\curvearrowright G$ and $\Lambda\curvearrowright H$ are orbit equivalent iff they are conjugate. Subsequently, this  has been generalized in \cite{Io14} to the case when $G$ and $H$ are arbitrary, not necessarily compact, connected Lie groups with trivial centers. The only difference is that, in the locally compact setting, the spectral gap assumption no longer makes sense and  has to be replaced with the assumption that the action $\Gamma\curvearrowright (G,{\bf m}_G)$ is {\it strongly ergodic}. After choosing a Borel probability measure $\mu$ on $G$ equivalent to ${\bf m}_G$, the latter requires that any sequence of measurable sets $A_n\subset G$ satisfying $\mu(gA_n\Delta A_n)\rightarrow 0$, for all $g\in\Gamma$, must be asymptotically trivial, in the sense that $\mu(A_n)\mu(A_n^c)\rightarrow 0$.

If $G$ is a connected compact simple Lie group (e.g., if $G=$ SO$(n)$, for $n\geq 3$),  the translation action $\Gamma\curvearrowright (G,{\bf m}_G)$ has spectral gap and thus is strongly ergodic, whenever $\Gamma$ is generated by matrices with algebraic entries. This result is due to Bourgain and Gamburd for $G=$ SO$(3)$ \cite{BG06}, and to  Benoist and  de Saxc\'{e} in general \cite{BdS14}. In joint work with Boutonnet and Salehi-Golsefidy, we have shown that the same holds if $G$ is an arbitrary connected simple Lie group \cite{BISG15}. In particular, the  translation action $\Gamma\curvearrowright (G,{\bf m}_G)$ is strongly ergodic, whenever $\Gamma<G:=$ SL$_n(\mathbb R)$ is a dense subgroup whose elements  are matrices with algebraic entries.

\section{Structure and  rigidity for group von Neumann algebras}
In this section, we survey some recent developments in the classification of von Neumann algebras arising from countable groups. Our presentation follows three directions.

\subsection{Structural results}\label{structure} We first present results which provide classes of group II$_1$ factors $L(\Gamma)$ with various indecomposability properties, such as primeness 
and the lack of Cartan subalgebras. 

In \cite{Po83},  Popa proved that II$_1$
factors arising from free groups with uncountably many generators are prime and do not have Cartan subalgebras. The first examples of  separable such II$_1$ factors were obtained
in the mid 1990s as an application of free probability theory.
Thus, Voiculescu showed that  the free group factors $L(\mathbb F_n)$, with $n\geq 2$, do not admit
Cartan subalgebras \cite{Vo95}. Subsequently, Ge used the techniques from \cite{Vo95} to prove that the free group factors are also prime \cite{Ge96}.

These results have been since generalized and strengthened in several ways. 
Using subtle C$^*$-algebras techniques, Ozawa remarkably proved that II$_1$ factors associated to  icc hyperbolic groups $\Gamma$ are {\it solid}: the relative commutant $A'\cap L(\Gamma)$ of any diffuse von Neumann subalgebra $A\subset L(\Gamma)$ is amenable \cite{Oz03}. In particular,  $L(\Gamma)$ and all of its non-amenable subfactors are prime.
By developing a novel technique based on closable derivations, Peterson showed that II$_1$ factors arising from icc  groups with positive first $\ell^2$-Betti number are prime \cite{Pe06}. A new proof of solidity of $L(\mathbb F_n)$ was found by Popa in \cite{Po06c} (see Section \ref{3}), while II$_1$ factors coming from icc groups $\Gamma$ admitting a proper cocycle into $\ell^2(\Gamma)$ were shown to be solid in \cite{Pe06}. 
For additional examples of prime and solid II$_1$ factors, see \cite{Oz04,Po06a,CI08,CH08,Va10b,Bo12,HV12,DI12,CKP14,Ho15}.

In \cite{OP07}, Ozawa and Popa discovered that the free group factors enjoy a  striking structural property, called {\it strong solidity}, which strengthens both primeness and absence of Cartan subalgebras: the normalizer of any diffuse amenable subalgebra $A\subset L(\mathbb F_n)$ is amenable. Generalizing this result, Chifan and Sinclair  showed that in fact the von Neumann algebra of any icc hyperbolic group is strongly solid \cite{CS11}.
For further examples of strongly solid factors, see  \cite{OP08,Ho09,HS09,Si10}. 

The groups $\Gamma$ for which $L(\Gamma)$ was shown to be prime have certain properties, such as hyperbolicity or the existence of unbounded quasi-cocycles, relating them to lattices in rank one Lie groups. On the other hand, the primeness question for  the higher rank  lattices such as PSL$_m(\mathbb Z)$, $m\geq 2$, remains a major open problem. Moreover, little is known about the structure of II$_1$ factors associated to lattices in higher rank semisimple Lie groups. We formulate a general question in this direction:

\vskip 0.05in
{\bf Problem V}\label{V}. {\it Let $\Gamma$ be an icc irreducible lattice in a direct product $G=G_1\times...\times G_n$ of connected non-compact simple real Lie groups with finite center. 
Prove that $L(\Gamma)$ is prime and does not have a Cartan subalgebra. Moreover, if $G_1,...,G_n$ are of rank one, prove that $L(\Gamma)$ is strongly solid.}
\vskip 0.05in

A subgroup $\Gamma<G$ is called a {\it lattice} if it is discrete and has finite co-volume, ${\bf m}_G(G/\Gamma)<\infty$. A lattice $\Gamma<G$ is called {\it irreducible} if its projection onto   $\times_{j\not=i}G_j$ is dense, for all $1\leq i\leq n$.

When $n=1$, Problem V has been resolved in \cite{Oz03} and \cite{CS11} (see also \cite{OP08,Si10}) whenever $G=G_1$  has rank one, but is open if $G$ has higher rank. 

In the rest of the section, we record recent progress in the case $n\geq 2$.
Note that the irreducibility assumption on $\Gamma$ is needed in order to exclude product lattices $\Gamma=\Gamma_1\times...\times\Gamma_n$,  whose II$_1$ factors are obviously non-prime.
However, it seems plausible that $L(\Gamma)$ does not have a Cartan subalgebra, for an arbitrary lattice $\Gamma$ of $G$, irreducible or not. 
Indeed, this was established by Popa and Vaes \cite{PV12} if $G_1,...,G_n$ are all of rank one. 
Complementing their work, the first examples of prime II$_1$ factors arising from lattices in higher rank semisimple Lie groups were only recently obtained in \cite{DHI16}:

\begin{theorem}\emph{\cite{DHI16}}
\label{DHI}
If $\Gamma$ is an icc irreducible lattice in a product $G = G_1\times...\times G_n$ of $n\geq 2$ connected non-compact rank one simple real Lie groups with finite center, then $L(\Gamma)$ is prime.
\end{theorem}

Theorem \ref{DHI} implies that the II$_1$ factor associated to PSL$_2(\mathbb Z[\sqrt{2}])$, which can be realized as an irreducible lattice in PSL$_2(\mathbb R)\times$PSL$_2(\mathbb R)$, is prime. 
Theorem \ref{DHI} also applies when each $G_i$ is a rank one non-compact simple algebraic group over a local field. This implies that the II$_1$ factor arising from PSL$_2(\mathbb Z[\frac{1}{p}])$, which is an irreducible lattice in PSL$_2(\mathbb R)\times$PSL$_2(\mathbb Q_p)$, is prime for  any prime $p$.
It remains however an open problem whether PSL$_2(\mathbb Z[\sqrt{2}])$,  PSL$_2(\mathbb Z[\frac{1}{p}])$, or any other group covered by Theorem \ref{DHI} gives rise to a strongly solid II$_1$ factor.

To overview the proof of Theorem \ref{DHI}, put $M=L(\Gamma)$ and assume the existence of a decomposition $M=M_1\bar{\otimes}M_2$ into a tensor product of  factors.
Let $\Delta_{\Gamma}:M\rightarrow M\bar{\otimes}M$ be the embedding given by $$\Delta_{\Gamma}(u_g)=u_g\otimes u_g, \;\;\;\text{for all $g\in\Gamma$.}$$ 
Since $G_i$ has rank one, it admits a non-elementary hyperbolic lattice $\Lambda_i<G_i$.
As a consequence, $\Gamma$  is measure equivalent (in the sense of Gromov) to a product $\Lambda=\Lambda_1\times...\times\Lambda_n$ of hyperbolic groups.
In combination with the relative strong solidity of hyperbolic groups \cite{PV12}, this allows one to conclude that a corner of $\Delta_{\Gamma}(M_j)$  embeds into $M\bar{\otimes}M_j$, for all $j\in\{1,2\}$.  
By making crucial use of an ultrapower technique from \cite{Io11a}, we derive
that $\Gamma$ admits commuting non-amenable subgroups. Indeed, by \cite{Io11a} the existence of commuting non-amenable subgroups can be deduced whenever we can find II$_1$ subfactors $A,B\subset M$ such that a corner of $\Delta_{\Gamma}(A)$ embeds into $M\bar{\otimes}(B'\cap M)$.
While more work is needed in general, this easily gives a contradiction  for  $\Gamma=\text{PSL}_2(\mathbb Z[\sqrt{2}])$ or PSL$_2(\mathbb Z[\frac{1}{p}])$.

\subsection{Algebraic rigidity}\label{algebraic}
The primeness results discussed in Section \ref{structure} provide classes of groups for which the absence of a direct product decomposition is inherited by their von Neumann algebras. In a complementary direction, it was recently  shown in \cite{CdSS15} that  for a wide family of product groups $\Gamma$,  the von Neumann algebra $L(\Gamma)$ completely remembers the product structure of $\Gamma$.
 
\begin{theorem}\emph{\cite{CdSS15}}\label{CdSS}
Let $\Gamma=\Gamma_1\times...\times\Gamma_n$, where  $\Gamma_1,...,\Gamma_n$ are $n\geq 2$ icc hyperbolic groups.
Then any countable group $\Lambda$ such that $L(\Gamma)\cong L(\Lambda)$ is a product of $n$ icc groups, $\Lambda=\Lambda_1\times...\times\Lambda_n$.
\end{theorem}

By applying the unique prime factorization theorem from \cite{OP03}, Theorem \ref{CdSS} can be strengthened to moreover show that $L(\Lambda_i)$ is stably isomorphic to $L(\Gamma_i)$, for all $i$. 

To outline the proof of Theorem \ref{CdSS}, assume that $n=2$ and put $M=L(\Gamma)$ and $M_i=L(\Gamma_i)$.
Consider the embedding $\Delta_{\Lambda}:M\rightarrow M\bar{\otimes}M$ associated to an arbitrary group von Neumann algebra decomposition $M=L(\Lambda)$. 
First, solidity properties of hyperbolic groups \cite{Oz03,BO08} are used to derive the existence of $i,j\in\{1,2\}$ such that a corner of $\Delta_{\Lambda}(M_i)$ embeds into $M\bar{\otimes}M_j$. The ultrapower technique of \cite{Io11a} then allows to transfer the presence of non-amenable commuting subgroups from $\Gamma$ to the mysterious group $\Lambda$.
Using a series of ingenious combinatorial lemmas and strong solidity results from \cite{CS11,PV12},  the authors conclude that $\Lambda$ is indeed a product of icc groups.

Motivated by Theorem \ref{CdSS}, it seems natural to investigate what other constructions in group theory can be detected at the level of the associated von Neumann algebra. 
Very recently,  progress on this problem has been made in \cite{CI17} by providing a class of amalgamated free product groups 
 $\Gamma=\Gamma_1*_{\Sigma}\Gamma_2$ whose von Neumann algebra $L(\Gamma)$  entirely recognizes the amalgam structure of $\Gamma$.
 This class includes all groups of the form $(A*B\times A*B)*_{A\times A}(A*B\times A*B)$, where $A$ is an icc amenable group and $B$ is a non-trivial hyperbolic group.
For such groups $\Gamma$ it is shown in \cite{CI17} that any  group $\Lambda$ satisfying $L(\Gamma)\cong L(\Lambda)$ admits a amalgamated free product  decomposition $\Lambda=\Lambda_1*_{\Omega}\Lambda_2$ such that the inclusions $L(\Sigma)\subset L(\Gamma_i)$ and $L(\Omega)\subset L(\Lambda_i)$ are isomorphic, for all $i\in\{1,2\}$.

\subsection{W$^*$-superrigidity} 
The  rigidity results presented in Sections \ref{structure} and \ref{algebraic} give instances when various algebraic aspects of groups can be recovered from their von Neumann algebras. 
We now discuss the most extreme type of rigidity for group von Neumann algebras. This occurs when the  von Neumann algebra $L(\Gamma)$ completely remembers the group $\Gamma$. Thus, we say that a countable group $\Gamma$ is {\it W$^*$-superrigid} if any group $\Lambda$ satisfying $L(\Gamma)\cong L(\Lambda)$ must be isomorphic to $\Gamma$. 

The first class of W$^*$-superrigid groups was discovered in our joint work with Popa and Vaes \cite{IPV10}:

\begin{theorem}\emph{\cite{IPV10}}\label{IPV}
Let $G_0$ be any non-amenable group and $S$ be any infinite amenable group. Define the wreath product group $G=G_0^{(S)}\rtimes S$, and consider the left multiplication action of $G$ on $I=G/S$.  Then the generalized wreath product group $\Gamma=(\mathbb Z/2\mathbb Z)^{(I)}\rtimes G$ is W$^*$-superrigid.
\end{theorem}

The conclusion of Theorem \ref{IPV} does not hold for {\it plain} wreath product groups $\Gamma=(\mathbb Z/2\mathbb Z)^{(G)}\rtimes G$. 
In fact, for any non-trivial torsion free group $G$, there is a torsion free group $\Lambda$ with $L(\Gamma)\cong L(\Lambda)$, while $\Gamma$ and $\Lambda$ are not isomorphic. 
Nevertheless, for a large family of groups $G$, including all icc property (T) groups, we are able to conclude that any group $\Lambda$ satisfying $L(\Gamma)\cong L(\Lambda)$ decomposes as a semi-direct product $\Gamma=\Sigma\rtimes G$, for some abelian group $\Sigma$.
This is in particular recovers a seminal result of Popa \cite{Po04a}  showing that if $\Lambda=(\mathbb Z/2\mathbb Z)^{(H)}\rtimes H$ is also a plain wreath product, then  $L(\Gamma)\cong L(\Lambda)$ entails $G\cong H$. 

To describe the strategy of the proof of Theorem \ref{IPV}, let $M=L(\Gamma)$ and assume that $M=L(\Lambda)$, for a countable group $\Lambda$. 
Then we have two embeddings $\Delta_{\Gamma},\Delta_{\Lambda}:M\rightarrow M\bar{\otimes}M$. 
 Notice that $M$ can also be realized as the group measure space factor of the generalized Bernoulli action $G\curvearrowright \{0,1\}^{I}$.
By extending the methods of \cite{Io10} from plain to generalized Bernoulli actions, we give a classification of all  possible embeddings $\Delta:M\rightarrow M\bar{\otimes}M$.  This enables us to deduce the existence of a unitary element $u\in M\bar{\otimes}M$ such that $\Delta_{\Lambda}(x)=u\Delta_{\Gamma}(x)u^*$, for all $x\in M$.   A principal novelty of \cite{IPV10} is being able to conclude that the groups $\Gamma$ and $\Lambda$ are isomorphic (in fact, unitarily conjugate modulo scalars inside $M$) from the mere existence of $u$.

Following \cite{IPV10}, several other classes of W$^*$-superrigid groups were found in \cite{BV12,CI17}.  Thus, it was shown in \cite{BV12} that the left-right wreath product $\Gamma=(\mathbb Z/2\mathbb Z)^{(\Gamma)}\rtimes(\Gamma\times\Gamma)$ is W$^*$-superrigid, for any icc hyperbolic group $\Gamma$. 
Very recently, a class of W$^*$-superrigid amalgamated free product groups was found in \cite{CI17}. These groups are also C$_{\text{red}}^*$-superrigid since, unlike the examples from \cite{IPV10,BV12}, they do not admit non-trivial normal amenable subgroups.

While the above results provide several large families of W$^*$-superrigid groups, the superrigidity question remains open for many natural classes of groups including the lattices PSL$_m(\mathbb Z)$, $m\geq 3$. In fact, a well-known {\it rigidity conjecture} of Connes \cite{Co82} asks if  $L(\Gamma)\cong L(\Lambda)$ for icc property (T) groups $\Gamma$ and $\Lambda$ implies $\Gamma\cong \Lambda$. Since property (T) is a von Neumann algebra invariant \cite{CJ83}, Connes' conjecture is equivalent to asking whether icc property (T) groups $\Gamma$ are W$^*$-superrigid.


\end{document}